\documentclass[11pt]{article}

\usepackage{amsthm,amsmath,amssymb}

{\theoremstyle{definition}
    \newtheorem{thm}{Theorem}[section]
    \newtheorem{prop}[thm]{Proposition}
    \newtheorem{lem}[thm]{Lemma}
    
    \newtheorem{defn}[thm]{Definition}
    \newtheorem{note}[thm]{Note}
    \newtheorem{alg}[thm]{Algorithm}
    \newtheorem{rem}[thm]{Remark}
    \newtheorem{exmp}[thm]{Example}
    
}

\newcommand{\OO}{{\mathcal {O}}}

\newcommand{\weak}{{\curlyvee}}

\DeclareMathOperator{\codim}{codim}

\DeclareMathOperator{\Hcodim}{Hcodim}

\DeclareMathOperator{\Ker}{Ker}
\DeclareMathOperator{\Sing}{Sing}
\DeclareMathOperator{\MaxB}{\underline{\mathbf{Max}}}

\DeclareMathOperator{\Eord}{E-ord}

\DeclareMathOperator{\Sat}{Sat}
\DeclareMathOperator{\Spec}{Spec}

\DeclareMathOperator{\inv}{inv}

\def\der#1#2{\dfrac{\partial #1}{\partial #2}}

\title{Embedded desingularization of toric varieties.}

\author{Roc\'{\i}o Blanco and
Santiago Encinas\footnote{This research was partially supported by MTM2006-10548, 
MTM2007-64704}}

\begin{document}

\maketitle

\begin{abstract}

We present a new method to achieve an embedded desingularization of a
toric variety.

Let \(W\) be a regular toric variety defined by a fan \(\Sigma\) and 
\( X\subset W \) be a toric embedding.
We construct a finite sequence of combinatorial blowing-ups such that 
the final strict transforms \( X'\subset W' \) are regular and \( X' 
\) has normal crossing with the exceptional divisor.

\end{abstract}



\section*{Introduction}

Fix a polynomial ring, a toric (not necessarily normal) variety is defined by a 
prime ideal generated by binomials. Such varieties can be considered 
as combinatorial, in fact all the information they carry can be 
expressed in terms of combinatorial objects. This gives a way of 
computing geometric invariants of the toric variety.
There are also many applications of the theory of toric varieties, 
see for example \cite{Cox1997}. For an introduction to toric varieties see 
\cite{Danilov1978}, \cite{Oda1988} or \cite{Fulton1993}.

This paper is devoted to construct an algorithm of desingularization 
of toric varieties and log-resolution of binomial ideals.
Given a binomial prime ideal, corresponding to a toric variety \( X \), we 
construct a sequence of combinatorial blowing-ups such that the 
strict transform of the variety \( X' \) is non singular and has 
normal crossings with the exceptional divisor.
Our algorithm is valid if the ground field is perfect of any 
characteristic.

In \cite{GonzalezTeissier2002} an algorithm of embedded desingularization of toric 
varieties is given. In this paper the authors construct a toric map 
\( X'\to X \), which provides an embedded desingularization of \( X 
\). 
If \( X\subset W \) is a toric embedding, they construct a refinement
of the fan of \( W \) depending on the embedding, see also
\cite{Teissier2004}.  By \cite{ConciniProcesi1983,ConciniProcesi1985} such a toric
morphism \( X'\to X \) may be dominated by a sequence of blowing-ups
along regular centers.

There is another desingularization method of toric varieties 
\cite{BierstoneMilman2006}, producing a sequence of blowing-ups along 
non singular centers.
This method is defined in terms of the Hilbert function of the 
variety.

In \cite{BlancoSomolinos2008thesis,BlancoSomolinos2009prep} an 
algorithm of log-resolution of binomial ideals is constructed based 
on the computation of an ordering function \( \Eord \).
The function \( \Eord \) is the order of the ideal defining the 
variety along a normal crossing divisor \( E \).
But this log-resolution algorithm
depends on a choice of a Gr\"obner basis of the ideal.
\medskip

The algorithm presented here depends on the ordering function \( \Eord
\) and a codimension function \( \Hcodim \).  This algorithm does not
depend on any choice and it can be implemented at the computer and we
expect to have a working implementation shortly.

On the other hand, we prove an equivalence of the geometric notion of 
tranversality of a variety with respect to a normal crossing divisor 
and a new notion of transversality of \( \mathbb{Z} \)-modules. 
Then we are able to translate geometric notions to combinatorial terms.

This can be considered as a first step on a more ambitious program 
translating notions from toric varieties, in terms of dual cones and 
fans, to notions in terms of the binomial equations of the variety.

We want to thank to Ignacio Ojeda for useful suggestions and 
conversations to improve the presentation of this paper. We are also 
grateful to Bernard Teissier and Pedro 
Gonz\'alez for useful comments.
\medskip

The paper is structured as follows:
First we recall known facts on toric varieties. There is a bijection 
between toric varieties and saturated \( \mathbb{Z} \)-submodules of 
\( \mathbb{Z}^{n} \). We define a notion of transversality of \( 
\mathbb{Z} \)-submodules which will be equivalent to the usual notion of 
transversality of the toric variety with respect to a normal crossing 
divisor. Given an affine toric variety \( X \), we will prove also 
the existence of a minimal regular toric variety \( V \) tranversal 
to \( E \) (the normal crossing divisor) and 
containing \( X \).
This minimal embedding \( X\subset V \) will define a function \( 
\Hcodim \) which is the first coordinate of our resolution function.
The rest of the coordinates of this resolution function comes from 
the process of \( E \)-resolution constructed in 
\cite{BlancoSomolinos2008thesis} and \cite{BlancoSomolinos2009prep}.
In the last section we will construct the algorithm of embedded 
desingularization of a toric variety \( X\subset W \), where \( W \) 
is a smooth toric variety.

\section*{Notation and first definitions}

\begin{note} \label{FirstDef}
Fix a perfect field \( k \).
We will denote the affine space of dimension \( n \) as usual \(
\mathbb{A}^{n}=\Spec(k[x_{1},\ldots,x_{n}]) \).
The torus of dimension \( n \) is \( 
\mathbb{T}^{n}=\Spec(k[x_{1}^{\pm},\ldots,x_{n}^{\pm}]) \).
Note that \( \mathbb{A}^{n}\setminus\mathbb{T}^{n} \) is a union of 
\( n \) hypersurfaces having only normal crossings.
Let \( E \) be a set of hypersurfaces such that every \( H\in E \) is 
an irreducible component of \( \mathbb{A}^{n}\setminus\mathbb{T}^{n} \).
The set \( E \) corresponds to a subset of indexes \( 
\mathbf{E}\subset\{1,\ldots,n\} \).
We set
\begin{equation*}
    \mathbb{TA}^{n}_{E}=\mathbb{A}^{n}\setminus
    \bigcup_{H\not\in E} H=
    \Spec\left(k[x_{1},\ldots,x_{n}]_{\prod_{i\not\in\mathbf{E}}x_{i}}\right)
\end{equation*}
where \( \Spec(k[x_{1},\ldots,x_{n}]_{\prod_{i\not\in\mathbf{E}}x_{i}}) \) is the 
localization of the polynomial ring with respect to the product \( 
\prod_{i\not\in E}x_{i} \).

Note that \( E \) is a set of hypersurfaces of
\( \mathbb{TA}^{n}_{E} \) having only normal crossings.
\end{note}

A morphism \( \mathbb{T}^{n}\to\mathbb{T}^{1} \) which is a group 
homomorphism is called a \emph{character}.

For example if \( a\in\mathbb{Z}^{n} \) then the morphism defined by 
\( T\to X_{1}^{a_{1}}\cdots X_{n}^{a_{n}} \) is a character.
In fact,
every character of \( \mathbb{T}^{n} \) is as above
\cite{CoxLittleSchenck2008prep,Humphreys1975,EisenbudSturmfels1996}.
So that the group of characters of a \( n \)-dimensional torus is a
free abelian group of rank \( n \).

\begin{note}
There is a bijection between morphisms \(
\mathbb{T}^{d}\to\mathbb{T}^{n} \) which are also group homomorphisms
and homomorphisms \( \mathbb{Z}^{n}\to\mathbb{Z}^{d} \) of \(
\mathbb{Z} \)-modules \cite{Humphreys1975}.  Moreover, closed reduced
immersions \( \mathbb{T}^{d}\to\mathbb{T}^{n} \) correspond to
surjective homomorphisms \( \mathbb{Z}^{n}\to\mathbb{Z}^{d} \).
\end{note}

\section{Affine toric varieties.}

We recall some basic definitions and well known results on toric varieties.

\begin{defn}
An affine toric variety is an affine variety \( X \) of dimension \( d
\), such that \( X \) contains the torus \( \mathbb{T}^{d} \) as a 
dense open set and the action of the torus extends to an action of \( 
\mathbb{T}^{d} \) to \( X \).
\end{defn}

Theorems~\ref{ThCaracAffTor} and \ref{ThRegAffTor} are well known 
results.
\begin{thm} \label{ThCaracAffTor}
\cite{CoxLittleSchenck2008prep,MillerSturmfels2005,EisenbudSturmfels1996}
Let \( X \) be a scheme of dimension \( d \).
TFAE:
\begin{enumerate}
    \item  \( X \) is an affine variety.

    \item  \( X\cong\Spec(k[t^{a_{1}},\ldots,t^{a_{n}}]) \), where \( 
    a_{1},\ldots,a_{n}\subset\mathbb{Z}^{d} \).

    \item  \( X\subset\mathbb{A}^{n} \) and \( I(X) \) is prime and 
    generated by binomials.
\end{enumerate}
\end{thm}

\begin{thm} \label{ThRegAffTor}
Let \( X \) be an affine toric variety of dimension \( d \).

\( X \) is regular if and only if \(
X\cong\mathbb{TA}^{d}_{E} \) for
some \( \mathbf{E}\subset\{1,\ldots,d\} \).
\end{thm}

\begin{defn} \label{DefAffTorEmbed}
An affine toric embedding is a reduced closed subscheme \( X\subset W \) where:
\begin{itemize}
    \item \( W \) is a regular affine toric variety.

    \item  \( X \) is an affine toric variety.

    \item The inclusion is toric, which means that one has a group
    homomorphism from the torus of \( X \) to the torus of \( W \).
\end{itemize}
\end{defn}

Toric varieties are related to \( \mathbb{Z} \)-submodules of \( 
\mathbb{Z}^{n} \).
\begin{defn} \label{DefSaturaM}
Let \( M\subset\mathbb{Z}^{n} \) be a \( \mathbb{Z} \)-module. The 
saturation of \( M \) is:
\begin{equation*}
    \Sat(M)=\{\alpha\in\mathbb{Z}^{n}\mid \lambda\alpha\in M
    \text{ for some } \lambda\in\mathbb{Z}\}
\end{equation*}
We say that a \( \mathbb{Z} \)-module \( L\subset\mathbb{Z}^{n} \) 
is  saturated if \( \Sat(L)=L \).
\end{defn}

Note that \( L\subset\mathbb{Z}^{n} \) is saturated if and only if the 
quotient \( \mathbb{Z}^{n}/L \) is a free \( \mathbb{Z} \)-module.

Note also that \( M\otimes_{\mathbb{Z}}\mathbb{Q}=
\Sat(M)\otimes_{\mathbb{Z}}\mathbb{Q} \).
\medskip

The following theorem is based on known results.
\begin{thm} \cite{EisenbudSturmfels1996}
Let be \( n\in\mathbb{N} \), \( r\leq n \) and \( d\leq n \).
There is a bijection correspondence between the following sets:
\begin{enumerate}
    \item  The set of affine toric embeddings \( 
    X\subset\mathbb{TA}^{n}_{E} \), with \( d=\dim{X} \).

    \item  The set of closed and reduced immersions \( 
    \mathbb{T}^{d}\to\mathbb{T}^{n} \).
    
    \item  The set of surjective homomorphisms of \( \mathbb{Z} 
    \)-modules, \( \mathbb{Z}^{n}\to\mathbb{Z}^{d} \).

    \item  The set of saturated \( \mathbb{Z} \)-submodules \( L\subset 
    \mathbb{Z}^{n} \) of rank \( n-d \).
\end{enumerate}
\end{thm}

\section{Sublattices of \( \mathbb{Z}^{n} \)}

In this section we introduce the notion of transversality of a \( 
\mathbb{Z} \)-module with respect to a subset \( 
\mathbf{E}\subset\{1,\ldots,n\} \).
We will prove that there exists always a maximal transversal submodule of 
any saturated \( \mathbb{Z} \)-submodule of \( \mathbb{Z}^{n} \).

\begin{lem} \label{LemaZCodim1}
Let \( \alpha=(\alpha_{1},\ldots,\alpha_{n})\in\mathbb{Z}^{n} \) with
\( \alpha\neq 0 \).

The following are equivalent:
\begin{enumerate}
    \item  \( \gcd\{\alpha_{1},\ldots,\alpha_{n}\}=1 \).

    \item The \( \mathbb{Z} \)-module \( (\alpha) \) is saturated
    (\ref{DefSaturaM}).
\end{enumerate}
Moreover assume that \( \alpha_{m+1}=\cdots=\alpha_{n}=0 \) for some
\( m\leq n \) and that \( (\alpha) \) is saturated, then there is a
surjective homomorphism \( \psi:\mathbb{Z}^{n}\to \mathbb{Z}^{n-1}
\) with \( \Ker\psi=(\alpha) \) and
such that \( \psi(e^{(n)}_{j})=e^{(n-1)}_{j-1} \) for \( 
j=m+1,\ldots,n \), where \( e^{(n)}_{j}\in\mathbb{Z}^{n} \) is the \( 
j \)-th element of the canonical base of \( \mathbb{Z}^{n} \).
\end{lem}

\begin{proof}
The equivalence of (1) and (2) is an easy exercise.

The last assertion follows from the fact that the Smith normal form 
of the column matrix \( \alpha \) is \( (1,0,\ldots,0) \).
There is a (non unique) invertible integer 
matrix \( P \) such that \( P\alpha=(1,0,\ldots,0) \).
In fact if \( \alpha_{m+1}=\cdots=\alpha_{n}=0 \) then \( P \) may be 
chosen as follows
\begin{equation*}
    P=\left(
    \begin{array}{cc}
        P' & 0  \\
        0 & I_{n-m}
    \end{array}\right)
\end{equation*}
Set \( A \) the \( (n-1)\times n \) matrix obtained by deleting the 
first row of \( P \). The matrix \( A \) defines the required 
homomorphism \( \mathbb{Z}^{n}\to\mathbb{Z}^{n-1} \).
\end{proof}

\begin{defn} \label{DefZE}
Let \( \mathbf{E}\subset\{1,\ldots,n\} \) be a subset.
We define \( \mathbb{Z}^{n}_{\mathbf{E}} \) and \( 
\mathbb{Z}^{n}_{\mathbf{E}^{+}} \) to be
\begin{equation*}
    \mathbb{Z}^{n}_{\mathbf{E}}=\{
    (\alpha_{1},\ldots,\alpha_{n})\in \mathbb{Z}^{n} \mid
    \alpha_{j}\geq 0\ \forall j\in \mathbf{E}
    \}
\end{equation*}
\begin{equation*}
    \mathbb{Z}^{n}_{\mathbf{E}^{+}}=\{
    (\alpha_{1},\ldots,\alpha_{n})\in \mathbb{Z}^{n} \mid
    \alpha_{j}> 0\ \forall j\in \mathbf{E}
    \}
\end{equation*}
\end{defn}

\begin{defn} \label{DefLatticeTrans}
Let \( M\subset\mathbb{Z}^{n} \) be a \( \mathbb{Z} \)-submodule and 
\( \mathbf{E}\subset\{1,\ldots,n\} \) be a subset.

We say that \( M \) is weak-transversal to \( \mathbf{E} \) if \( M \) admits 
a system of \( \mathbb{Z} \)-generators \( 
\alpha_{1},\ldots,\alpha_{\ell} \) with \( \alpha_{j}\in 
\mathbb{Z}^{n}_{\mathbf{E}} \) (\ref{DefZE}). 

We say that a \( \mathbb{Z} \)-module \( L\subset\mathbb{Z}^{n} \) is
transversal to \( \mathbf{E} \) if it is weak-transversal to \(
\mathbf{E} \) and \( L \) is saturated.
\end{defn}

\begin{defn} \label{DefEM}
Let \( M\subset\mathbb{Z}^{n} \) be a \( \mathbb{Z} \)-submodule and 
\( \mathbf{E}\subset\{1,\ldots,n\} \) be a subset.
Consider \( pr_{j}:\mathbb{Z}^{n}\to\mathbb{Z} \) the \( j 
\)-projection, \( j=1,\ldots,n \).
Set
\begin{equation*}
    \mathbf{E}_{M}=\{j\in\mathbf{E}\mid pr_{j}(M)\neq 0\}
\end{equation*}
\end{defn}

\begin{rem} \label{RemNodeg}
The set \( \mathbf{E}_{M} \) depends on the module \( M \), but we can
reduce the study of transversality of \( M \) to the smaller subset \(
\mathbf{E}_{M} \).  Note that for any generator system \(
\alpha_{1},\ldots,\alpha_{\ell} \) of \( M \) we have that \(
\alpha_{i,j}=0 \) for all \( i=1,\ldots,\ell \) and any \( j\in
\mathbf{E}\setminus \mathbf{E}_{M} \).
\end{rem}

Propositions~\ref{PropReduceEM} and \ref{PropReduceGamma} come 
from discussions with Ignacio Ojeda.
\begin{prop} \label{PropReduceEM}
Let \( M\subset\mathbb{Z}^{n} \) be a \( \mathbb{Z} \)-submodule and 
\( \mathbf{E}\subset\{1,\ldots,n\} \) be a subset.

The \( \mathbb{Z} \)-module \( M \) is weak transversal to \(
\mathbf{E} \) (\ref{DefLatticeTrans}) if and only if \( M \) is weak
transversal to \( \mathbf{E}_{M} \).
\end{prop}

\begin{proof}
One implication is obvious.
Assume that \( M \) is weak transversal to \( \mathbf{E}_{M} \).
So that there exists a generator system \( 
\alpha_{1},\ldots,\alpha_{\ell} \) such that \( \alpha_{i,j}\geq 0 \) 
for \( i=1,\ldots,\ell \) and \( j\in\mathbf{E}_{M} \).
The result follows from remark~\ref{RemNodeg}.
\end{proof}

\begin{prop} \label{PropReduceGamma}
Let \( M\subset\mathbb{Z}^{n} \) be a \( \mathbb{Z} \)-submodule and 
\( \mathbf{E}\subset\{1,\ldots,n\} \) be a subset.

The \( \mathbb{Z} \)-module \( M \) is weak transversal to \( 
\mathbf{E} \) if and only if there is \( \gamma\in 
M\cap\mathbb{Z}^{n}_{\mathbf{E}_{M}^{+}} \)
\end{prop}

\begin{proof}
Assume that \( M \) is weak transversal to \( E \). There is a 
generator system \( \alpha_{1},\ldots,\alpha_{\ell} \) of \( M \) with \( 
\alpha_{i}\in\mathbb{Z}^{n}_{\mathbf{E}} \), \( i=1,\ldots,\ell \).

Note that for any \( j\in\mathbf{E}_{M} \) there is an index \( 
i\in\{1,\ldots,\ell\} \) such that \( \alpha_{i,j}>0 \).
Set \( \gamma=\alpha_{1}+\cdots+\alpha_{\ell} \) and it is clear that 
\( \gamma\in\mathbb{Z}^{n}_{\mathbf{E}_{M}^{+}} \).

Conversely, assume that there is
\( \gamma\in M\cap\mathbb{Z}^{n}_{\mathbf{E}_{M}^{+}} \).
Consider \( \beta_{1},\ldots,\beta_{\ell} \) a generator system of \( 
M \).
There are integers \( a_{1},\ldots,a_{\ell}\in\mathbb{Z} \) with
\( \gamma=a_{1}\beta_{1}+\cdots+a_{\ell}\beta_{\ell} \). We may 
assume that \( \gcd\{a_{1},\ldots,a_{\ell}\}=1 \).
Note that we may complete \( \gamma \) to a generator system of \( M 
\), say, \( \gamma,\gamma_{2},\ldots,\gamma_{\ell} \).
This is a consequence of the fact that the Smith normal form of the 
row matrix \( (a_{1},\ldots,a_{\ell}) \) is \( (1,0,\ldots,0) \).

Now we may choose positive integers \( 
\lambda_{2},\ldots,\lambda_{\ell} \) such that \( 
\gamma_{i}+\lambda_{i}\gamma\in\mathbb{Z}^{n}_{\mathbf{E}} \).
Set \( \alpha_{1}=\gamma \) and \( 
\alpha_{i}=\gamma_{i}+\lambda_{i}\gamma \), \( i=2,\ldots,\ell \).
It is clear that \( \alpha_{1},\ldots,\alpha_{\ell} \) is a generator 
system of \( M \) and \( \alpha_{i}\in\mathbb{Z}^{n}_{\mathbf{E}} \), 
\( i=1,\ldots,\ell \). 
In fact we may assume that \( 
\alpha_{i}\in\mathbb{Z}^{n}_{\mathbf{E}_{M}^{+}} \).
\end{proof}

\begin{prop} \label{PropZModuleL}
Let \( \mathbf{E}\subset\{1,\ldots,n\} \) and
let \( M\subset\mathbb{Z}^{n} \) be a \( \mathbb{Z} \)-module.

If \( M \) is weak-transversal to \( \mathbf{E} \) then \( \Sat(M) \) 
is transversal to \( \mathbf{E} \).
\end{prop}

\begin{proof}
Note that \( \mathbf{E}_{M}=\mathbf{E}_{\Sat(M)} \).
Then proposition~\ref{PropZModuleL}
is a direct consequence of \ref{PropReduceEM} and
\ref{PropReduceGamma}.
\end{proof}

\begin{prop} \label{ExistLatTrans}
Let \( \mathbf{E}\subset\{1,\ldots,n\} \) be a subset.
Let \( L\subset \mathbb{Z}^{n} \) be a saturated \( \mathbb{Z} \)-submodule.

There exists a unique \( \mathbb{Z} \)-module \( L_{0} \) such that 
\begin{itemize}
    \item  \( L_{0}\subset L \),

    \item  \( L_{0} \) is transversal to \( \mathbf{E} \),

    \item  If \( L'_{0}\subset L \) and \( L'_{0} \) is transversal 
    to \( \mathbf{E} \) then \( L'_{0}\subset L_{0} \).
\end{itemize}
\end{prop}

\begin{proof}
Consider all \( \mathbb{Z} \)-submodules \( 
\{M_{\lambda}\}_{\lambda\in\Lambda} \) such that \( 
M_{\lambda}\subset L \) and \( M_{\lambda} \) is weak-transversal to 
\( \mathbf{E} \) for every \( \lambda\in\Lambda \).

Set \( M=\sum_{\lambda\in\Lambda}M_{\lambda} \). Note that \( M \) 
is weak-transversal to \( \mathbf{E} \) and \( M\subset L \).

By \ref{PropZModuleL} \( L_{0}=\Sat(M) \) is
transversal to \( \mathbf{E} \) and we have also that \( M\subset L_{0}\subset L \).
In fact by construction \( M=L_{0} \) and it is the biggest \( 
\mathbb{Z} \)-module with this property.
\end{proof}

\section{Affine toric varieties and transversality.}

In this section we will prove the equivalence of the new  notion of 
transversality of \( \mathbb{Z} \)-submodules 
(definition~\ref{DefLatticeTrans}) and the geometric usual notion of 
transversality of a variety with respect to a normal crossing divisor.
\medskip

Let \( W \) be a regular affine toric variety of dimension \( n \).
It follows from (\ref{ThRegAffTor}) that \( W\cong \mathbb{TA}^{n}_{E} \)
(notation as in \ref{FirstDef}).
Recall that \( E \) is a set of regular hypersurfaces in \( W \) 
having only normal crossings. Using the isomorphism \( W\cong 
\mathbb{TA}^{n}_{E} \) we may identify \( E \) with a set \( 
\mathbf{E}\subset\{1,\ldots,n\} \).
With this identification
\begin{equation*}
    W=\mathbb{TA}^{n}_{E}=
    \Spec(k[x_{1},\ldots,x_{n}]_{\prod_{i\not\in\mathbf{E}}x_{i}}).
\end{equation*}

The variety \( W=\mathbb{TA}^{n}_{E} \) has a distinguished point \(
\xi_{0}\in W \)
\begin{equation*}
    \xi_{0}\in\bigcap_{H\in E}H
\end{equation*}
with coordinates \( \xi_{0}=(\xi_{0,1},\ldots,\xi_{0,n}) \) where \(
\xi_{0,i}=0 \) if \( i\in \mathbf{E} \) and \( \xi_{0,i}=1 \) if \(
i\not\in\mathbf{E} \).

\begin{defn} \label{DefEord}
    \cite{BlancoSomolinos2008thesis,BlancoSomolinos2009prep}
Let \( W \) be a regular affine toric variety of dimension \( n \) 
and let \( \mathcal{J}\subset\OO_{W} \) be a sheaf of ideals.
For any point \( \xi\in W \) consider \( E_{\xi} \) the intersection of all 
hypersurfaces \( H\in E \) with \( \xi\in H \):
\begin{equation*}
    E_{\xi}=\bigcap_{\xi\in H\in E}H
\end{equation*}
The ideal \( I(E_{\xi})\subset\OO_{W} \) is generated by all the 
equations of hypersurfaces \( H \) with \( \xi\in H\in E \).

We define the function \( \Eord(J):W\to\mathbb{N} \) as follows:
\begin{equation*}
    \Eord(J)(\xi)=\max\{b\in\mathbb{N}\mid J\subset I(E_{\xi})^{b}\}
\end{equation*}
where \(J\) is an ideal in \(W\).
\end{defn}

Note that the function \( \Eord(J) \) is constant along the strata
defined by \( E \).  In fact \( \Eord(J)(\xi) \) is the (usual) order
of the ideal \( J \) at the generic point of \( E_{\xi} \).  The
function \( \Eord(J):W\to\mathbb{N} \) is upper-semi-continuous, see
\cite{BlancoSomolinos2008thesis,BlancoSomolinos2009prep} for a proof
and more details.

\begin{defn}
An ideal \( J\subset\OO_{W} \), \( W=\mathbb{TA}^{n}_{E} \), is 
binomial if \( J \) can be generated, as ideal, by binomials: \( 
x^{\alpha}-x^{\beta} \), with \(\alpha,\beta\in \mathbb{N}^n\).
\end{defn}

\begin{lem} \label{LemaDistinguido}
Let \( J\subset\OO_{W} \) be a 
binomial ideal, \( W=\mathbb{TA}^{n}_{E} \), and let \( \xi_{0}\in W \) 
be the distinguished point.

Then \( \xi_{0}\in\MaxB\Eord(J)=\{\xi\in W\mid \Eord(J)(\xi)=\max\Eord(J)\} \).
\end{lem}

\begin{proof}
Note that \( E_{\xi_{0}}\subset E_{\xi} \) for any \( \xi\in W \).
\end{proof}

The following definition is general for any variety: 
\begin{defn} \label{DefTransVar}
Let \( X\subset W \) be an embedded variety and let \( E \) be a set of 
regular hypersurfaces of \( W \) having only normal crossings.

We say that \( X \) is transversal to \( E \) at a point \( \xi\in X 
\) if there is a regular system of parameters of \( \OO_{W,\xi} \), 
\( x_{1},\ldots,x_{n}\in\OO_{W,\xi} \), such that
\begin{itemize}
    \item  \( I(X)_{\xi}=(x_{1},\ldots,x_{r}) \) for some \( r\leq n 
    \) and

    \item For all \( H\in E \) with \( \xi\in H \), then \(
    I(H)_{\xi}=(x_{i}) \) for some \( i \) with \( r<i\leq n \).
\end{itemize}
\end{defn}

Consider \( W=\mathbb{TA}^{n}_{E} \), for some \(
\mathbf{E}\subset\{1,\ldots,n\} \).
The derivatives with poles along \( E \) is a free \( \OO_{W} \)-module
of rank \( n \) and a natural basis of this module is
\begin{equation*}
    x_{i}^{\epsilon_{i}}\der{\ }{x_{i}}\qquad i=1,\ldots,n
\end{equation*}
where \( \epsilon_{i}=0 \) if \( i\in\mathbf{E} \) and \(
\epsilon_{i}=1 \) if \( i\not\in\mathbf{E} \).

\begin{lem} \label{LemaJacMatrix}
Let \( X\subset W=\mathbb{TA}^{n}_{E} \) be an affine toric embedding 
with \( d=\dim(X) \).

Consider any set of binomial generators of the ideal \( 
I(X)\subset\OO_{W} \)
\begin{equation*}
    I(X)=(x^{\alpha_{1}^{+}}-x^{\alpha_{1}^{-}},\ldots, 
    x^{\alpha_{m}^{+}}-x^{\alpha_{m}^{-}}) =(f_{1},\ldots,f_{m})
\end{equation*}
Fix a point \( \xi\in X \).
The variety \( X \) is transversal to \( E \) at the point \( \xi \)
(\ref{DefTransVar}) if and only if the 
jacobian matrix:
\begin{equation*}
    \left(x_{i}^{\epsilon_{i}}\der{f_{j}}{x_{i}}\right)_{i,j}
\end{equation*}
has rank \( n-d \) at the point \( \xi \).
\end{lem}

\begin{proof}
This lemma is a direct consequence of a general fact on algebraic
varieties.
\end{proof}

\begin{prop} \label{PropEordPos}
Let \( X\subset W=\mathbb{TA}^{n}_{E} \) an affine toric embedding.

If \( \max\Eord(I(X))>0 \) then \( X \) is not transversal to \( E \).
\end{prop}

\begin{proof}
It follows from \ref{LemaJacMatrix} and \ref{LemaDistinguido}.
At the distinguished point \( \xi_{0} \),
the jacobian matrix in lemma~\ref{LemaJacMatrix} is zero modulo the maximal ideal
at \( \xi_{0} \).
\end{proof}

\begin{thm} \label{ThVTrans1}
Let \( V\subset W=\mathbb{TA}^{n}_{E} \) be an affine toric embedding.

They are equivalent:
\begin{enumerate}
    \item  \( V \) is transversal to \( E \).

    \item The ideal \( I(V) \) is generated by hyperbolic equations
    \begin{equation*}
	I(V)=(x^{\alpha_{1}}-1,\ldots,x^{\alpha_{\ell}}-1)
    \end{equation*}
    where \( \ell=n-\dim{V} \), \(
    \alpha_{1},\ldots,\alpha_{\ell}\in\mathbb{Z}^{n}_{\mathbf{E}} \) and they
    generate a saturated lattice of rank \( \ell \).
\end{enumerate}
\end{thm}

\begin{proof}
Set \( \ell=n-\dim{V} \).

Let \( \alpha_{1},\ldots,\alpha_{\ell}\in
\mathbb{Z}^{n}_{\mathbf{E}} \) be such that \( 
\alpha_{1},\ldots,\alpha_{\ell} \) they generate a saturated lattice
of rank \( \ell \).
Assume that \( I(V)=(x^{\alpha_{1}}-1,\ldots,x^{\alpha_{\ell}}-1) \).
Consider the jacobian matrix (\ref{LemaJacMatrix})
\begin{equation*}
    \left(x_{i}\frac{\partial\ }{x_{i}}\left(x^{\alpha_{j}}-1\right)\right)=
    \left(
    \begin{array}{ccc}
        \alpha_{1,1}x^{\alpha_{1}} & \cdots &
	\alpha_{\ell,1}x^{\alpha_{\ell}}  \\
        \vdots &  & \vdots  \\
        \alpha_{1,n}x^{\alpha_{1}} & \cdots &
	\alpha_{\ell,n}x^{\alpha_{\ell}}
    \end{array}
    \right)
\end{equation*}
Note that the rank of this matrix at any point \( \xi\in V \) is the 
rank of the matrix \( \left(\alpha_{1}|\cdots |\alpha_{\ell}\right) 
\) having \( \alpha_{i} \) as columns.
And the rank of this matrix is \( \ell \) (independently of the 
characteristic of the ground field \( k \)).
So that \( V \) is transversal to \( E \).
\medskip

Conversely assume that \( V \) is transversal to \( E \). We may
assume that \( \mathbf{E}=\{r+1,\ldots,n\} \).

Let show first the codimension one case: \( \ell=1 \).
So that \( I(V)=(x^{\alpha^{+}_{1}}-x^{\alpha^{-}_{1}}) \) where \( 
\alpha_{1}=\alpha^{+}_{1}-\alpha^{-}_{1} \) and \( 
\alpha^{+}_{1},\alpha^{-}_{1}\in\mathbb{Z}^{r}\times\mathbb{N}^{n-r} \).
By \ref{PropEordPos} we have that \( \max\Eord(I(V))=0 \) so that we 
may assume that \( I(V)=(x^{\alpha_{1}}-1) \) with \( 
\alpha_{1}\in\mathbb{Z}^{r}\times\mathbb{N}^{n-r} \).
Since \( V \) is a toric variety the ideal \( I(V) \) is prime and \( 
\gcd\{\alpha_{1,1},\ldots,\alpha_{1,n}\}=1 \). Lemma~\ref{LemaZCodim1} 
gives the result.

We now prove the general case, codimension \( \ell>1 \).  By
\ref{PropEordPos} there is an hyperbolic equation \( x^{\alpha_{1}}-1\in
I(V) \).  Since \( V \) is toric, we may assume that \(
\gcd\{\alpha_{1,1},\ldots,\alpha_{1,n}\}=1 \).  Set \( W_{1} \) the toric
hypersurface defined by \( x^{\alpha_{1}}-1 \).  After reordering the last
\( n-r \) coordinates we may assume that \(
\alpha_{1,m+1}=\cdots=\alpha_{1,n}=0 \) and \( \alpha_{1,i}>0 \) for \(
i=r+1,\ldots,m \), where \( r\leq m \).  We do not assume anything on \(
\alpha_{1,1},\ldots,\alpha_{1,r} \).

Let \(
\psi:\mathbb{Z}^{n}\to\mathbb{Z}^{n-1} \) be the homomorphism given by
lemma~\ref{LemaZCodim1}.  Note that \( W_{1}\cong
\mathbb{T}^{m-1}\times\mathbb{A}^{n-m} \).
We have \( V\subset W_{1} \) and by induction there are \( 
\bar{\beta}_{2},\ldots,\bar{\beta}_{\ell}\in 
\mathbb{Z}^{m-1}\times\mathbb{N}^{n-m} \) such that the ideal of \( V 
\) in \( W_{1} \) is generated by \( 
y^{\bar{\beta}_{2}}-1,\ldots,y^{\bar{\beta}_{\ell}}-1 \), and \( 
\bar{\beta}_{2},\ldots,\bar{\beta}_{\ell} \) generate a saturated
lattice in \( \mathbb{Z}^{n-1} \) of rank \( \ell-1 \).
Let \( \beta_{2},\ldots,\beta_{\ell}\in\mathbb{Z}^{n} \) such that \( 
\psi(\beta_{i})=\bar{\beta}_{i} \) for \( i=2,\ldots,\ell \).
We have that \( \alpha_{1},\beta_{2},\ldots,\beta_{\ell} \) generate a
saturated lattice of rank \( \ell \).

It is clear that \( \beta_{i}\in\mathbb{Z}^{m}\times\mathbb{N}^{n-m} 
\), but in general \( 
\beta_{i}\not\in\mathbb{Z}^{r}\times\mathbb{N}^{n-r} \).
Since \( \alpha_{1,i}>0 \) for \( i=r+1,\ldots,m \) there are natural 
numbers \( \lambda_{2},\ldots,\lambda_{\ell} \) such that \( 
\alpha_{i}=\beta_{i}+\lambda_{i}\alpha_{1}\in 
\mathbb{Z}^{r}\times\mathbb{N}^{n-r} \), \( i=2,\ldots,\ell \).
And \( \alpha_{1},\alpha_{2},\ldots,\alpha_{\ell} \) generate the same
lattice.
Finally we have the equality of ideals
\begin{equation*}
    I(V)=(x^{\alpha_{1}}-1,x^{\alpha_{2}}-1,\ldots,x^{\alpha_{\ell}}-1)
\end{equation*}
\end{proof}

\begin{lem} \label{LemaHiperIdeal}
Let \( J=(x^{\beta_{1}}-1,\ldots,x^{\beta_{s}}-1) \) be an ideal 
generated by hyperbolic binomials, \( \beta_{i}\in 
\mathbb{Z}^{n}_{\mathbf{E}} \), \( i=1,\ldots,s \).
Assume that \( J \) is a prime ideal.

If \( \gamma\in\mathbb{Z}^{n} \)
then \( x^{\gamma^{+}}-x^{\gamma^{-}}\in J \) if and only if \( 
\gamma \) belongs to the \( \mathbb{Z} \)-module generated by
\( \beta_{1},\ldots,\beta_{s} \) in \(\mathbb{Z}^{n} \).
\end{lem}

\begin{proof}
The \( \mathbb{Z} \)-module \( L \) generated by \(
(\beta_{1},\ldots,\beta_{s}) \) is associated to the toric variety
defined by \( J \).

It is well known that \( \beta\in L \) if and only if \( 
x^{\beta^{+}}-x^{\beta^{-}}\in J \).
\end{proof}

\begin{prop} \label{EquivTransLatV}
Let \( V\subset\mathbb{TA}^{n}_{E} \) be an affine toric embedding.
Set \( L\subset\mathbb{Z}^{n} \) be the lattice associated to \( V \).

\( V \) is transversal to \( E \) if and only if \( L \) is 
transversal to \( \mathbf{E} \).
\end{prop}

\begin{proof}
It is a consequence of \ref{ThVTrans1} and \ref{LemaHiperIdeal}.
\end{proof}

\begin{thm}\label{ThExistVMin}
Let \( X\subset W=\mathbb{TA}^{n}_{E} \) be an affine toric 
embedding.

There is a unique toric variety \( V \) such that the embeddings \( 
X\subset V\subset W \) are toric and \( V \) is the minimum toric 
variety containing \( X \) and transversal to \( E \).
\end{thm}

\begin{proof}
It follows from \ref{ExistLatTrans} and \ref{EquivTransLatV}.
\end{proof}

\section{Embedded toric varieties.}

In the previous sections we have reduced to the case of affine toric 
varieties. We generalize here to (non affine) toric varieties, and
we define the first coordinate of our resolution function.
\medskip

Let \( W \) be a regular toric variety defined by a fan \( \Sigma \) 
\cite{Fulton1993}.
Let \( T\subset W \) be the torus of \( W \), which is open and dense 
in \( W \).
Set \( E \) the simple normal crossing divisor given by \( W\setminus 
T \).

For every \( \sigma\in\Sigma \) the open set \( W_{\sigma}\subset W 
\) is an affine toric variety, so that \(
W_{\sigma}\cong\mathbb{TA}^{n}_{E} \).

\begin{defn} \label{DefToricEmbed}
A toric embedding is a closed subscheme \( X\subset W \) such that 
for every \( \sigma\in\Sigma \) if \( X_{\sigma}=X\cap W_{\sigma} 
\) then \( X_{\sigma}\subset W_{\sigma} \) is an affine toric 
embedding (\ref{DefAffTorEmbed}).
\end{defn}

For every \( \sigma\in\Sigma \) there is a unique toric affine 
variety \( V_{\sigma}\subset W_{\sigma} \) transversal to \( E \) and 
such that \( X_{\sigma}\subset V_{\sigma} \) (theorem~\ref{ThExistVMin}).
\medskip

In fact, this toric affine variety \( V_{\sigma} \) is a
regular toric affine variety.

\begin{rem} \label{RemSigmaXi}
Note that for any \( \xi\in W \), there is a unique \( \sigma\in\Sigma
\) such that \( \xi\in W_{\sigma} \) and the affine open set \(
W_{\sigma} \) is minimum with this property.  In fact \( \xi \)
belongs to the orbit of the distinguished point of \( W_{\sigma} \).
\end{rem}

\begin{defn} \label{DefHCodim}
Let \( \xi\in X \subset W\) be a point.  Let \( X_{\sigma} \) be the
minimum affine open set containing the point \( \xi \).

The hyperbolic codimension of \( X \) at \( \xi \) is
\begin{equation*}
    \Hcodim(X)(\xi)=\dim{V_{\sigma}}-\dim{X}
\end{equation*} 
where \( V_{\sigma}\subset W_{\sigma} \) is the minimum toric affine
variety such that \( V_{\sigma}\supset X_{\sigma} \) and it is
transversal to \( E \) (\ref{ThExistVMin}).
\end{defn}

\begin{rem}
The hyperbolic codimension \( \Hcodim(X) \) (\ref{DefHCodim}) can be
understood as a toric embe\-dding dimension.  The number \( \Hcodim(X)
\) at \( \xi \) is the minimum dimension of a regular toric variety \(
V \) including \( X \).

In the case \(V_{\sigma}=W_{\sigma}\), then
\(\Hcodim(X)(\xi)=\codim_{W}(X)\), the codimension of \(X\) in \(W\).
\end{rem}

\begin{rem} \label{CombCenterW}
\cite{BierstoneMilman2006}
Let \( \Delta\in\Sigma \) be an element of the fan \( \Sigma \)
defining the regular variety \( W \).
The cone \( \Delta \) defines a smooth closed subvariety \(
Z_{\Delta}\subset W \) as follows:

The toric variety \( W \) is covered by affine toric varieties \(
W_{\sigma} \) with \( \sigma\in\Sigma \).  So that \( Z_{\Delta} \) is
covered by affine pieces \( (Z_{\Delta})_{\sigma}=Z_{\Delta}\cap
W_{\sigma} \), \( \sigma\in\Sigma \).

If \( \Delta \) is not a face of \( \sigma \) then \(
(Z_{\Delta})_{\sigma}=\emptyset \).

If \( \Delta \) is a face of \( \sigma \), note that \(
W_{\Delta}\subset W_{\sigma} \) is an open inclusion.
Then \( (Z_{\Delta})_{\sigma} \) is the (closure) of the orbit of the 
distinguished point of \( W_{\Delta} \).

The smooth closed center \( Z_{\Delta} \) we will say that it is a 
\emph{combinatorial center} of \( W \). In fact note that at every 
affine chart \( W_{\sigma}\cong \mathbb{TA}^{n}_{E} \), for some \( E 
\), the combinatorial center \( Z_{\Delta}\cap W_{\sigma} \) is 
defined by some coordinates \( x_{i} \) with \( i\in\mathbf{E} \).
\end{rem}

\begin{rem}
Note that if \(X\subset W\cong \mathbb{TA}^n_{E}\) is an affine toric embedding and 
\(Z_{\Delta}\subset W\) is a combinatorial center, then the strict 
transforms \(X'\subset W'\) give an affine toric embedding.    
\end{rem}

\begin{prop} \label{HcodimDecrease}
Let \( \Delta\in\Sigma \) and \( Z_{\Delta}\subset W \) the 
combinatorial center associated to \( \Delta \) (\ref{CombCenterW}).

Let \( W'\to W \) be the blow-up with center \( Z_{\Delta} \). Set \( 
X'\subset W' \) the strict transform of \( X \).
If \( \xi'\in X' \) then
\begin{equation*}
    \Hcodim(X')(\xi')\leq \Hcodim(X)(\xi)
\end{equation*}
where \( \xi' \) maps to \( \xi \).
\end{prop}

\begin{proof}
Let \(W_{\sigma}\) be the minimum affine open set of \( W \) containing the point
\( \xi \) and let \( V_{\sigma} \) be the minimum toric affine 
variety in \( W_{\sigma} \)
such that \( V_{\sigma}\supset X_{\sigma} \) and it is transversal to
\( E \) (\ref{ThExistVMin}).

Let \( W'_{\sigma'} \), with \( \sigma'\in\Sigma' \), be the minimum
affine open set of \( W' \) containing the point \( \xi' \).
Let \( V'_{\sigma'} \) be the minimum toric affine variety such
that \( V'_{\sigma'}\supset X'_{\sigma'} \) and it is transversal to
\( E' \).
Note that \( X'_{\sigma'}\subset(X_{\sigma})' \) is an open immersion,
where \( (X_{\sigma})'\subset X' \) is the strict transform of \(
X_{\sigma}\subset X \).

Let \( \left(V_{\sigma}\right)' \) be the strict transform of \( V_{\sigma} \).
Note that \( \left(V_{\sigma}\right)' \) is smooth and transversal to \( E' \).
So that \( \left(V_{\sigma}\right)'\cap W'_{\sigma'}\supset 
V'_{\sigma'} \).
And the result follows from the last inclusion.
\end{proof}

\section{\( E \)-resolution of binomial ideals.}

In \cite{BlancoSomolinos2008thesis,BlancoSomolinos2009prep} were given
some notions in terms of \emph{binomial basic objects along E}, where
\(E\) was a normal crossing divisor in the ambient space \(W\).
In terms of the \( \Eord \) (\ref{DefEord}) one may construct a 
sequence of combinatorial blowing-ups such that the transform of a 
given binomial ideal has maximal \( E \)-order equal to zero.

We remind here the main results.  For more details on the several
constructions and proofs see \cite{BlancoSomolinos2008thesis} or
\cite{BlancoSomolinos2009prep}.
All these notions work for general binomial ideals,
without any restriction.

\begin{note}
Using this structure of binomial basic object along \(E\), and the
language of \emph{mobiles} (see \cite{EncinasHauser2002}), it is
possible to construct a resolution function involving the \(E\)-order
of certain ideals computed by induction on the dimension of \(W\).
\end{note}

\begin{note} \label{induction}
Roughly speaking, given \((W,(J,c),H,E)\), where \(J\) is a binomial
ideal, \(c\) is a positive integer, and \(H\) is the set of
exceptional hypersurfaces, by induction on the dimension of \(W\),
construct ideals \(J_i\) defined in local flags \(W=W_n\supseteq
W_{n-1}\supseteq \cdots \supseteq W_i\supseteq \cdots \supseteq W_1\),
and then objects \((W_i,(J_i,c_{i+1}),H_i,E_i)\) in dimension \(i\),
where each \(E_i=W_i\cap E\), \(H_i=W_i\cap H\).  The integer numbers
\(c_{i+1}\) are computed as the \(E\)-order of certain ideals
\(P_{i+1}\) coming from the previous dimension \(i+1\), that is
\(c_{i+1}=\max\ \Eord(P_{i+1})\) is the \emph{critical value} in
dimension \(i\).  Denote \(c_{n+1}=c\).  \\

If the \emph{\(E\)-singular locus} of \((J_i,c_{i+1})\) is non empty,
then factorize the ideal \(J_i=M_i\cdot I_i\), where each ideal
\(M_i\) is defined by a normal crossings divisor supported by the
current exceptional locus \(H_i\).
\end{note}

\begin{note}
To make this induction on the dimension of \(W\), in
\cite{BlancoSomolinos2008thesis,BlancoSomolinos2009prep} it is proved
the existence of hypersurfaces of \emph{\(E\)-maximal contact} at any
stage of the resolution process.  This hypersurfaces are always
coordinate hyperplanes, and produce a combinatorial center to be blown
up.  Combinatorial centers are convenient to preserve the binomial
structure of the ideal after blow-up.
\end{note}

\begin{defn} \label{bboe}
A \emph{binomial basic object along \(E\)} is a tuple
\(B=(W,(\mathcal{J},c),H,E)\) where
\begin{itemize}
	\item \(W\) is a regular toric variety defined by a fan \( \Sigma \).
	
	\item \(E\) is the simple normal crossing divisor given by \( 
	W\setminus T \), where \( T\subset W \) is the torus of \( W \).
	
	\item \((\mathcal{J},c)\) is a \emph{binomial pair}, this
	means that \(\mathcal{J}\subset \mathcal{O}_W\) is a coherent
	sheaf of binomial ideals with respect to \(E\), and \(c\) is a
	positive integer number.  Note that for any \( \sigma\in\Sigma
	\) the sheaf of ideals \( \mathcal{J} \) restricted to the
	open affine subset \( W_{\sigma}\subset W \) is a binomial
	ideal \( J\neq 0 \) in \(
	k[x_{1},\ldots,x_{n}]_{\prod_{i\not\in\mathbf{E}}x_{i}} \).
		
	\item \(H\subset E\) is a set
	of normal crossing regular hypersurfaces in \(W\).
\end{itemize}
\end{defn}  

\begin{defn}
Let \(J\subset \mathcal{O}_{W}\) be a binomial ideal, \(c\) a positive
integer.  We call \emph{\(E\)-singular locus} of \(J\) with respect to
\(c\) to the set, $$E\text{-}\Sing(J,c)=\{\xi \in W/\
\Eord_{\xi}(J)\geq c\}.$$
\end{defn}

\begin{rem}
The \(E\)-singular locus is a closed subset of \(W\).
\end{rem}

\begin{defn}
Let \(J\subset \mathcal{O}_W\) be a binomial ideal.  Let \(\xi\in W\)
be a point such that \(\Eord_{\xi}(J)=\max\ \Eord(J)=\theta_E\).  A
hypersurface \(V\) is said to be a hypersurface of
\emph{\(E\)\text{-}maximal contact} for \(J\) at the point \(\xi\) if
\begin{itemize}
	\item [-] \(V\) is a regular hypersurface, \(\xi\in V\),

	\item [-] \(E\text{-}\Sing(J,\theta_E) \subseteq V\) and their
	transforms under blowing up with a combinatorial center \(Z_{\Delta}\subset
	V\) also satis\-fy \(E\text{-}\Sing(J',\theta_E) \subseteq
	V'\), whereas the \(E\)-order, \(\theta_E\) remains constant.
	\\ That is, \(\Eord_{\xi}(J')=\Eord_{\xi}(J)\), where \(J'\)
	is the \emph{controlled} transform of \(J\) and \(V'\) is the strict
	transform of \(V\). 
\end{itemize}
\end{defn}  

\begin{rem}
The controlled transform of \(J\) is the ideal
\(J'=I(Y')^{-\theta_E}\cdot J^{*}\) where \(Y'\) is the exceptional
divisor and \(J^{*}\) is the total transform of \(J\) under blowing
up.
\end{rem}

\begin{prop}
Let \(J\subset \mathcal{O}_W\) be a binomial ideal.  There exists a
hypersurface of \(E\)\text{-}maximal contact for \(J\).
\end{prop}

\begin{defn} \label{Einv} 
Let \((W,(J,c),H,E)\) be a binomial basic object along \(E\).  For all
point \(\xi\in E\text{-}\Sing(J,c)\) the resolution function
\(E\text{-}\inv_{(J,c)}\) will have \(n\) components with
lexicographical order, and will be of one of the following types:
$$E\text{-}\inv_{(J,c)}(\xi)=
\left\{\begin{array}{ll}
\left(\frac{\Eord_{\xi}(I_n)}{c_{n+1}},
\frac{\Eord_{\xi}(I_{n-1})}{c_n},\ldots,
\frac{\Eord_{\xi}(I_{n-r})}{c_{n-r+1}},\
\infty,\ \infty,\ldots,\infty\right) & (a) \vspace*{0.2cm} \\
\left(\frac{\Eord_{\xi}(I_n)}{c_{n+1}},
\frac{\Eord_{\xi}(I_{n-1})}{c_n},\ldots,
\frac{\Eord_{\xi}(I_{n-r})}{c_{n-r+1}},
\Gamma(\xi),\infty,\ldots,\infty\right) & (b)  \vspace*{0.2cm} \\
\left(\frac{\Eord_{\xi}(I_n)}{c_{n+1}},
\frac{\Eord_{\xi}(I_{n-1})}{c_n},\ldots,
\frac{\Eord_{\xi}(I_{n-r})}{c_{n-r+1}},\ldots,
\frac{\Eord_{\xi}(I_{1})}{c_2}\right) & (c)
\end{array}\right. $$ 

where the ideals \(I_i\) and the integer numbers \(c_i\) are as in
note \ref{induction}.  \\

In the case \(J_i=1\), for some \(i< n\), define
\(({E\text{-}\inv_{(J,c)}}_{i}(\xi),\ldots,
{E\text{-}\inv_{(J,c)}}_1(\xi))=(\infty,\ldots,\infty)\) in order to
preserve the number of components.

If \(\Eord_{\xi}(I_i)=0\), for some \(i< n\), then
\({E\text{-}\inv_{(J,c)}}_i(\xi)=\Gamma(\xi)\), where \(\Gamma\) is
the resolution function corresponding to the \emph{monomial case}, see
\cite{EncinasVillamayor2000, BravoEncinasVillamayor2005}.
And complete the resolution function
with the needed number of \(\infty\) components.
\end{defn}

\begin{note}
The \(E\text{-}\inv_{(J,c)}\) function is an upper-semi-continuous
function, see \cite{BlancoSomolinos2008thesis,BlancoSomolinos2009prep}.
We also denote \(
E\text{-}\inv_{(J,c)}(\xi)=E\text{-}\inv_{\xi}(J,c)\).

As a consequence of the upper-semi-continuity of the
\(E\text{-}\inv_{(J,c)}\) function,
$$\MaxB(E\text{-}\inv_{(J,c)})=\{\xi\in E\text{-}\Sing(J,c)|\
E\text{-}\inv_{(J,c)}(\xi)=\max\ E\text{-}\inv_{(J,c)}\}$$ is a closed
set.  In fact, it is the center of the next blow-up.

It can be proven that the \(E\text{-}\inv_{(J,c)}\) function drops
lexicographically after blow-up,
\cite{BlancoSomolinos2008thesis,BlancoSomolinos2009prep}.
\end{note}

\begin{lem} \label{dropEinv}
Let \( (W,(J,c),H,E)\) be a binomial basic object along \(E\).  Let
\(W \stackrel{\pi}{\leftarrow}W'\) be a blow-up with combinatorial center
\(Z_{\Delta}=\MaxB(E\text{-}\inv_{(J,c)}) \) then
$$E\text{-}\inv_{(J,c)}(\xi)>E\text{-}\inv_{(J',c)}(\xi')$$ where
\(\xi\in Z_{\Delta}\), \(\xi'\in Y'=\pi^{-1}(Z_{\Delta})\),
\(\pi(\xi')=\xi\).

The function \(E\text{-}\inv_{(J,c)}\) is the resolution function
associated to the binomial basic object along \(E\) given by
\((W,(J,c),H,E)\), and \(E\text{-}\inv_{(J',c)}\) corresponds to its
transform by the blow-up \(\pi\), \((W',(J',c),H',E')\).
\end{lem}

\begin{proof} See \cite{BlancoSomolinos2008thesis,BlancoSomolinos2009prep}.
\end{proof}

\begin{rem}
The \(E\text{-}\inv_{(J,c)}\) function provides a
\emph{\(E\)-resolution} of the binomial basic object along \(E\),
\((W,(J,c),H,E)\).
\end{rem}

\begin{defn}
Let \((W,(J,c),H,E)\) be a binomial basic object along \(E\), where 
\(E=\{E_{1},\ldots,E_{r}\}\), with \( r\leq n=\dim{W} \).
Let \(H=\{H_1,\ldots,H_s\}\subset E\) be the set of exceptional divisors, for
some \(s\leq r\).

We define a \emph{transformation} of the binomial basic object
$$(W,(J,c),H,E)\leftarrow (W',(J',c),H',E')$$ by means of the blowing
up \(W\stackrel{\pi}{\leftarrow} W'\), in a center
\(Z\subset E\text{-}\Sing(J,c)\), with
\begin{itemize}
    \item \(H'=\{H_1^{\curlyvee},\ldots,H_s^{\curlyvee},Y'\}\) where
    \(H_i^{\curlyvee}\), \( i=1,\ldots,s \), is the strict transform of
    \(H_i\) and \(Y'\)
    is the exceptional divisor in \(W'\).
    
    \item
    \(E'=\{E_{1}^{\curlyvee},\ldots,E_{r}^{\curlyvee},Y'\}\) where
    \(E_{i}^{\curlyvee}\), \( i=1,\ldots,r \), is the strict transform
    of \(E_{i}\) and \(Y'\) is the exceptional divisor in \(W'\).
    
    \item
    \(J'=I(Y')^{-c}\cdot J^{\ast}\) is the
    controlled transform of \(J\), where \(J^{\ast}\) is the
    total transform of \(J\).
\end{itemize}
\end{defn}

\begin{defn}
A sequence of transformations of binomial basic objects {\small
\begin{equation}
(W^{(0)}\!,(J^{(0)},c),H^{(0)}\!,E^{(0)})\!\leftarrow
\!(W^{(1)}\!,(J^{(1)},c),H^{(1)}\!,E^{(1)})\!\leftarrow\!\cdots
\leftarrow\!(W^{(N)}\!,(J^{(N)},c),H^{(N)}\!,E^{(N)})
\label{Eresolution} \end{equation}} is a \emph{\(E\)-resolution} of
\((W^{(0)},(J^{(0)},c),H^{(0)},E^{(0)})\), or simply a
\emph{\(E\)-resolution} of the pair \((J^{(0)},c)\), if
\(E\)-\(\Sing(J^{(N)},c)=\emptyset\).
\end{defn}

\begin{note}
This \(E\)-resolution function, the \(E\text{-}\inv_{(J,c)}\), works
for general binomial ideals, without any restriction, for more details
see \cite{BlancoSomolinos2008thesis,BlancoSomolinos2009prep}.

The \( E \)-resolution constructed in this way is independent of the
choice of coordinates.

In \cite{BlancoSomolinos2008thesis,BlancoSomolinos2009prep} was proved
that one may use this \( E \)-resolution in order to construct an
algorithm of log-resolution of binomial ideals and embedded
desingularization of binomial varieties.  But this algorithm of
log-resolution depends on a choice of a Gr\"obner basis of the
original binomial ideal.

In the next section we will construct an
algorithm (\ref{AlgoritmoTorico}) of embedded desingularization of toric
varieties which is
independent of the choice of coordinates.
\end{note}

\begin{thm} \label{ThEResol}
Let \(J\) be a binomial ideal.  If \(E\)-\(\Sing(J,c)\neq\emptyset\)
then there exists a \(E\)-resolution of \((J,c)\).
\end{thm}

\begin{proof}
The \(E\)-resolution of \((J,c)\) is given by the
\(E\text{-}\inv_{(J,c)}\) function, such that
\(E\)-\(\Sing(J^{\weak},c)=\emptyset\), this means
\(\max{\Eord(J^{\weak})}=0 \).
\end{proof}

\begin{rem} 
At any stage of the \(E\)-resolution process, the
\(E\text{-}\inv_{(J,c)}\) function determines the next combinatorial
center to be blown-up
\(Z_{\Delta}=\MaxB(E\text{-}\inv_{(J,c)})\), or equivalently,
the cone \(\Delta\), (\ref{CombCenterW}).
\end{rem}

\section{Embedded desingularization.}

Now we construct an algorithm of embedded desingularization of toric 
varieties. 
This algorithm is defined in terms of the function \( \Hcodim \) and 
a \( E \)-resolution of a suitable ideal, depending on the function 
\( E\text{-}\inv \).

\begin{prop} \label{drop}
Let \( X\subset W \) be a toric embedding and let \( \xi\in X \) be a
point.  Let \(X_{\sigma}\) be the minimum affine open set containing
the point \( \xi \) and let \( V_{\sigma} \) be the minimum toric
affine variety such that \( V_{\sigma}\supset X_{\sigma} \) and it is
transversal to \( E \).

At any affine open subset, \( X_{\sigma}\subset V_{\sigma}\subset
W_{\sigma}\).  Let \(W_{\sigma}
\stackrel{\pi}{\leftarrow}W'_{\sigma}\) be a blow-up with center
\(Z_{\Delta}\) then
$$(\Hcodim(X)(\xi),E\text{-}\inv_{\xi}(I_{V_{\sigma}}(X_{\sigma}),c))>
(\Hcodim(X')(\xi'),E\text{-}\inv_{\xi'}(I_{V'_{\sigma'}}(X'_{\sigma'}),c))$$

where \(\xi\in Z_{\Delta}\), \(\xi'\in Y'=\pi^{-1}(Z_{\Delta})\),
\(\pi(\xi')=\xi\), \(I_{V_{\sigma}}(X_{\sigma}) \) is the ideal of \(
X \) in \( V \), and \(c=\max \Eord(I_{V_{\sigma}}(X_{\sigma}))\).
\end{prop}

\begin{proof}
It follows from proposition (\ref{HcodimDecrease}) and
lemma~(\ref{dropEinv}).
\end{proof}
 
\begin{alg} \label{AlgoritmoTorico}
Let \( X\subset W \) be a toric embedding and let \( \xi\in X \) be a
point. Let \(W_{\sigma}\) be the minimum affine open set containing
the point \( \xi \) (\ref{RemSigmaXi}).
\begin{enumerate}
	\item Compute \( V_{\sigma} \), the minimum toric
affine variety such that \( V_{\sigma}\supset X_{\sigma} \) and it is
transversal to \( E \) (\ref{ThExistVMin}). 
	\item Set \(\Hcodim(X)(\xi)=\dim{V_{\sigma}}-\dim{X}\). 
\begin{itemize}
	\item If \(\Hcodim(X)(\xi)>0\) then compute
	\(E\text{-}\inv_{\xi}(I_{V_{\sigma}}(X_{\sigma}),1)\), here we
	set \(c=1\).  This determines \(Z_{\Delta}\).  Go to step (3).

	\item If \(\Hcodim(X)(\xi)=0\) the algorithm stops.  Note that
	in this case \(V_{\sigma}=X_{\sigma}\subset X\).
\end{itemize}
	\item Perform the blow-up with center \(Z_{\Delta}\) and go to step (1).
\end{enumerate}
\end{alg} 

Correctness of the algorithm follows by construction (\ref{DefHCodim} 
and \ref{HcodimDecrease}). Termination of
the algorithm follows from (\ref{drop}) and (\ref{ThEResol}).

\begin{rem}
The algorithm of embedded desingularization given in
\cite{BlancoSomolinos2008thesis} depends on the choice of a system of
coordinates.  Note that the new algorithm (\ref{AlgoritmoTorico})
given here does not depend on the coordinates election.
\end{rem}

\begin{thm} \label{ThEmbedDesing}
{\bf Embedded desingularization}

Let \( X\subset W \) be a toric embedding (\ref{DefToricEmbed}).  Let
\(E\) be the simple normal crossing divisor given by \(W\setminus T
\), where \( T\subset W \) is the torus of \( W \).

There exists a sequence of transformations of pairs
$$(W,E)\leftarrow (W^{(1)},E^{(1)})\leftarrow \cdots
\leftarrow (W^{(N)},E^{(N)})$$ which induces a proper birational
morphism \(\Pi: W^{(N)} \rightarrow W\) such that
\begin{enumerate}
	\item The restriction of this morphism \(\Pi\) to the
	\emph{regular locus of \(X\) along \(E\)}, defines an isomorphism
	$$Reg_E(X)\cong \Pi^{-1}(Reg_E(X))\subset W^{(N)}$$ where
	\(Reg_E(X)=\{\xi\in X|\ X \text{ is regular at } \xi \text{
	and has normal crossings with } E \}\).

	\item \(X^{(N)}\), the strict transform of \(X\) in
	\(W^{(N)}\), is regular and has normal crossings with the
	exceptional divisors \(E^{(N)}\).
\end{enumerate}
\end{thm}

\begin{proof}
It follows from correctness and termination of algorithm~\ref{AlgoritmoTorico}.
\end{proof}

The embedded desingularization of theorem~\ref{ThEmbedDesing} can be
implemented, since the key points are to determine the stratum
\(E_{\xi}\) (or the open subset \(W_{\sigma}\)) where the hyperbolic
codimension is maximum, and then to compute a combinatorial blow-up,
that can be easily encoded in the computer.

\begin{rem}
Theorem (\ref{ThEmbedDesing}) may be used to achieve a log-resolution 
of a toric ideal.

With the notation of theorem~(\ref{ThEmbedDesing}), let \( W^{(N+1)}\to 
W^{(N)} \) be the blowing up with center \( X^{(N)} \), which is a 
permissible center.
Note that the total transform of \( I(X)\OO_{W^{(N+1)}} \) is locally 
a monomial ideal (generated by monomials) and we may use an algorithm 
of log-resolution of monomial ideals as in \cite{Goward2005} or 
\cite{BierstoneMilman2006}.
So that log-resolution of the ideal \( I(X) \) follows from 
theorem~(\ref{ThEmbedDesing}) and log-resolution of monomial ideals.
\end{rem}

\begin{exmp}
Let \(X\subset W=\mathbb{A}^4\) be a toric embedding and let \( \xi\in X
\) be a point.  Let \(E=\{V(x),V(y),V(z),V(w)\}\) be the simple normal
crossing divisor given by \(W\setminus T \).  The toric variety \(X\)
is given by the equations $$X=\{x^2-y^3=0\}\cap\{xyz-w^2=0\}.$$ The
singular locus of the surface \(X\) is the \(z\)-axis.  Compute the
hyperbolic codimension at some points of \(X\), for example:
\begin{itemize}
	\item If \(\xi\not\in E\), then in a neighborhood of \(\xi\),
	\(W_{\sigma}=Spec(k[x^{\pm},y^{\pm},z^{\pm},w^{\pm}])\) and
	\(E_{\xi}=\emptyset\).  The minimum toric affine variety
	\(V_{\sigma}\) such that \( V_{\sigma}\supset X_{\sigma} \)
	and it is transversal to \( E_{\xi}\) is \(
	V_{\sigma}=X_{\sigma} \).  Then \(\Hcodim(X)(\xi)=0\).

	\item Set \( (V(x))^{c}=W\setminus V(x) \) be the complement of \(V(x)\).
	If \(\xi\in (V(x))^c\cap (V(y))^c\cap V(z)\cap V(w)\),
	then \(W_{\sigma}=Spec(k[x^{\pm},y^{\pm},z,w])\) and
	\(E_{\xi}=V(z)\cap V(w)\).  It is clear that \(
	V_{\sigma}=\{x^2-y^3=0\}\) and therefore
	\(\Hcodim(X)(\xi)=3-2=1\).
	
	\item If \(\xi\neq 0\) is a point at the \(z\)-axis,
	\(W_{\sigma}=Spec(k[x,y,z^{\pm},w])\) and \(E_{\xi}=V(x)\cap
	V(y)\cap V(w)\).  Assume \(\xi\) is the distinguished point of
	\(W_{\sigma}\), in this case \( V_{\sigma}=W_{\sigma}\) and
	\(\Hcodim(X)(\xi)=4-2=2\).
	
	\item If \(\xi\) is the origin,
	\(W_{\sigma}=Spec(k[x,y,z,w])\) and \(E_{\xi}=V(x)\cap
	V(y)\cap V(z)\cap V(w)\).  The minimum toric affine variety \(
	V_{\sigma}=W_{\sigma}\) and \(\Hcodim(X)(\xi)=4-2=2\).
\end{itemize}
It is easy to check that the hyperbolic codimension \( \Hcodim \) attains its
highest value along the \( z \)-axis.

If one computes the whole resolution function \(
(\Hcodim,E\text{-}\inv) \) (here we set \( c=1 \)) then its maximum value is
\begin{equation*}
    \max(\Hcodim(X),E\text{-}\inv)=
    (\Hcodim(X)(0),E\text{-}\inv(I_{V_{\sigma}}(X_{\sigma}),1))=
    (2,2,1,3/2,2)
\end{equation*}
and it is achieved at the origin, which is the first center to be
blown-up.
\medskip

We denote as \emph{\(x\)-th chart} the chart where we divide by \(x\).
For simplicity, we will denote each
\(\frac{y}{x},\frac{z}{x},\frac{w}{x}\) again as \(y,z,w\).  \medskip

At the \(x\)-th chart, the controlled transform of the ideal \(I(X)\) is  
\begin{equation*}
    I(X)'=x^{-1}\cdot(x^2-x^3y^3,x^3yz-x^2w^2)=x\cdot(1-xy^3,xyz-w^2).
\end{equation*}
If \(\eta'\) is a point that maps to the origin then \( \eta'\not\in X' 
\) and it lies in the first exceptional divisor \( V(x) \).
Consider \( \xi'\in X'\cap V(z)\cap V(w) \), the affine chart \( 
W'_{\sigma'} \) associated to \(\xi'\) as in  
(\ref{RemSigmaXi}) is
\(W'_{\sigma'}=Spec(k[x^{\pm},y^{\pm},z,w])\) and \(E'_{\xi'}=V(z)\cap
V(w)\).  The minimum toric affine variety \(V'_{\sigma'}\) containing
\(X'\) is \(V'_{\sigma'}=\{xy^3-1=0\}\) and
\(\Hcodim(X')(\xi')=3-2=1\).  The maximum value of the resolution
function is
\begin{equation*}
    \max(\Hcodim(X'),E\text{-}\inv)=
    (\Hcodim(X')(\xi'),E\text{-}\inv_{\xi'}(I_{V'_{\sigma'}}(X'_{\sigma'}),1))=
    (1,1,2,\infty,\infty)
\end{equation*}	
and it is reached along \(Z'=\{z=0\}\cap\{w=0\}\cap\{xy^3-1=0\}\).  Inside
\(V'_{\sigma'}\), the center is given by coordinates,
\(Z'_{\Delta}=\{z=0\}\cap\{w=0\}\), which is the next combinatorial center to be
blown-up.
After the blow-up at \(Z'_{\Delta}\), we consider the \(w\)-th chart 
\begin{equation*}
    I(X)''=w^{-1}\cdot(xyzw-w^2)=(xyz-w) \mod I(V'_{\sigma'}),
\end{equation*}
this means \(I(X)''=(xy^3-1,xyz-w)\).
\medskip

Let \( \xi''\in X'' \) mapping to \( \xi' \).
At this stage of the resolution process, the maximum value of the
resolution function is
\begin{equation*}
    \max(\Hcodim(X''),E\text{-}\inv)=
    (\Hcodim(X'')(\xi''),E\text{-}\inv_{\xi''}(I_{V''_{\sigma''}}(X''_{\sigma''}),1))=
    (1,1,1,\infty,\infty)
\end{equation*}	
and it is reached along
\(Z''=\{z=0\}\cap\{w=0\}\cap\{xy^3-1=0\}\).
After the blowing-up at \( Z''_{\Delta} \) we obtain two charts and 
for both \( \max\Hcodim(X''')=0 \). And \( X''' \) is regular and 
transversal to \( E''' \).
\end{exmp}

\bibliographystyle{amsalpha}

\providecommand{\bysame}{\leavevmode\hbox to3em{\hrulefill}\thinspace}

\end{document}